\documentclass[a4paper]{article}

\usepackage{amsmath}
\usepackage{amsthm}
\usepackage{amsfonts}
\newlength\PullBackLength
\newcommand\PullBack[1][2]{%
  \setlength{\global\PullBackLength}{#1em}%
  \kern\PullBackLength%
  &
  \kern-\PullBackLength}
\makeatletter
\@addtoreset{equation}{section}
\makeatother

\newtheorem{theorem}{Theorem}
\newtheorem{proposition}{Proposition}
\newtheorem{corollary}{Corollary}
\newtheorem{lemma}{Lemma}

\newcommand{\reals}{\mathbb{R}}
\newcommand{\vect}[1]{\mathbf{#1}}
\newcommand{\E}[1]{\mathbf{E}\left[#1\right]}

\DeclareMathOperator{\ran}{ran}
\DeclareMathOperator{\ind}{in}
\DeclareMathOperator{\out}{out}
\title{The Clustering Coefficient of a Scale-Free Random Graph}
\author{Nicole Eggemann\thanks{Supported by the EC Marie Curie programme NET-ACE (MEST-CT-2004-6724).}\mbox{ } and Steven D. Noble\thanks{Partially supported by the Heilbronn Institute for Mathematical Research, Bristol, UK.}\\Department of Mathematical Sciences, Brunel University\\Kingston Lane\\Uxbridge\\UB8 3PH\\United Kingdom\\ nicole.eggemann@brunel.ac.uk and steven.noble@brunel.ac.uk}
\begin{document}
\maketitle 
\begin{abstract}
We consider a random graph process in which, at each time step, a
new vertex is added with $m$ out-neighbours, chosen with
probabilities proportional to their degree plus a strictly
positive constant. We show that the expectation of the clustering
coefficient of the graph process is asymptotically proportional to
$\frac{\log n}{n}$. Bollob\'as and Riordan \cite{Bollobas02} have
previously shown that when the constant is zero, the same
expectation is asymptotically proportional to $\frac{(\log
n)^2}{n}$.
\end{abstract}

\section{Introduction}
Recently there has been a great deal of interest in the structure
of real world networks, especially the internet. Many mathematical
models have been proposed: most of these describe graph processes
in which new edges are added by some form of preferential
attachment. There is a vast literature discussing empirical
properties of these networks but there is also a growing body of
more rigorous work. A wide-ranging account of empirical properties
of networks can be found in \cite{Albert02}; a good survey of
rigorous results can be found in \cite{Bollobas02} or in the recent book \cite{Durret06}.

In \cite{Watts98} Watts and Strogatz defined `small-world'
networks to be those having small path length and being highly
clustered, and discovered that many real world networks are small-world
networks, e.g. the power grid of the western USA and the
collaboration graph of film actors.

There are conflicting definitions of the clustering coefficient appearing in the literature. See~\cite{Bollobas02} for a discussion of the relationships between them.
We define the clustering coefficient, $C(G)$ of a graph $G$ as follows:
\[ C(G) = \frac{3 \times \text{ number of triangles in }G}{\sum_{v\in V(G)} \binom{d(v)}{2}},\]
where $d(v)$ is the degree of vertex $v$.

The reason for the three in the numerator is to ensure that the clustering coefficient of a complete graph is one. This is the maximum possible value for a simple graph. However our graphs will not be restricted to simple graphs and so the clustering coefficient can exceed one. For instance if we take three vertices and join each pair by $m$ edges then the clustering coefficient is $m^2/(2m-1)$. Note that the clustering coefficient of a graph with at most $m$ edges joining any pair of vertices is at most $m$.

In this paper we establish rigorous results describing the
asymptotic behaviour of the clustering coefficient for one class of
model. Our graph theoretic
notation is standard. Since our graphs are growing, we let
$d_t(v)$ denote the total degree of vertex $v$ at time $t$. Sometimes we omit $t$ when the context is clear.

The Barab\'asi--Albert model (BA model)~\cite{Barabasi99} is perhaps the most widely studied graph process governed by preferential attachment. A new vertex is added to the graph at each time-step and is joined to $m$ existing vertices of the graph chosen with probabilities proportional to their degrees. A key observation~\cite{Barabasi99} is that in many large real-world networks, the proportion of vertices with degree $d$ obeys a power law.

In~\cite{Bollobas01} Bollob\'as et al. gave a mathematically
precise description of the BA model and showed rigorously that for
$d \le n^{\frac{1}{15}}$, the proportion of vertices with degree
$d$ asymptotically almost surely obeys a power law.

A natural generalisation of the BA model is to take the probability of attachment to $v$ at time $t+1$ to be proportional to $d_t(v)+a$, where $a$ is a constant representing the inherent attractiveness of a vertex. Buckley and Osthus~\cite{Buckley04} generalised the results in~\cite{Bollobas01} to the case where the attractiveness is a positive integer. A much more general model was introduced in~\cite{Cooper03} and further results extending~\cite{Bollobas01} were obtained.
Many more results on these variations of the basic preferential model can be found in~\cite{Bollobas02}.

Bollob\'as and Riordan showed~\cite{Bollobas02} that the
expectation of the clustering coefficient of the model
from~\cite{Bollobas01} is asymptotically proportional to $(\log
n)^2/n$. Bollob\'as and Riordan also considered
in~\cite{Bollobas02} a slight variant of the model
from~\cite{Bollobas01}. Their results imply that for this model
the expectation of the clustering coefficient is also
asymptotically proportional to $(\log n)^2/n$. We work with a
model depending on two parameters $\beta,m$, which to the best of
our knowledge was first studied rigorously by M\'ori in
\cite{Mori03}. In a sense, that we make precise in the next
section, Bollob\'as and Riordan's model is almost the special case
of M\'ori's model corresponding to $\beta=0$.

Our main result is to show that for $\beta>0$, asymptotically the
expectation of the clustering coefficient is proportional to $\log
n/n$. The main strategy of our proof follows~\cite{Bollobas02} and
we use very similar notation. In Section~\ref{sec:model} we give a
definition of the model that we use and explain its relationship
with the model studied in~\cite{Bollobas02}.
Section~\ref{sec:smallsub} contains results that give the
probability of the appearance of a small subgraph. We obtain the
expectation of the number of triangles appearing and of $\sum_v
\binom{d(v)}2$ in Section~\ref{sec:tri}. These two sections
follow~\cite{Bollobas02} quite closely. The overall aim is to
express the expectation of the clustering coefficient as the
quotient of the expectation of the number of triangles and the
expectation of $\sum_v \binom{d(v)}2$. We justify doing this in
Section~\ref{sec:final} and make use of a concentration result
proved in Section~\ref{sec:conc} using martingale methods.
Bollob\'as and Riordan~\cite{Bollobas02} used a similar strategy
and mentioned that they also used martingale methods.

\section{The model of M\'ori}\label{sec:model}
We now describe in detail M\'ori's generalisation of the BA model
\cite{Mori05}. Our definition involves a finer probability space
than was described in \cite{Mori05} but the underlying graph
process $(G^n_{m,\beta})$ is identical. The process depends on two
parameters: $m$ the outdegree of each vertex except the first and
$\beta \in \reals$ such that $\beta > 0$. (In \cite{Mori05},
M\'ori imposed the weaker condition that $\beta > -1$).

We first define the process when $m=1$. Let $G^1_{1,\beta}$
consist of a single vertex $v_1$ with no edges. The graph
$G^{n+1}_{1,\beta}$ is formed from $G^n_{1,\beta}$ by adding a new
vertex $v_{n+1}$ together with a single directed edge $e$. The
tail of $e$ is $v_{n+1}$ and the head is determined by a random
variable $f_{n+1}$. We diverge slightly from \cite{Mori05} in our
description of $f_{n+1}$.

Label the edges of $G^n_{1,\beta}$ with $e_2,\ldots,e_n$ so that
$e_i$ is the unique edge whose tail is $v_i$. Now let \[\Omega_{n+1} =
\{(1,v),\ldots,(n,v),(2,h),\ldots,(n,h),(2,t),\ldots,(n,t)\}.\] We
define $f_{n+1}$ to take values in $\Omega_{n+1}$ so that for $1 \leq i
\leq n$,
\[\Pr(f_{n+1} = (i,v)) = \frac{\beta}{(2+\beta)n-2}\]
and for $2 \leq i \leq n$,
\[  \Pr(f_{n+1} = (i,h)) = \Pr(f_{n+1}=(i,t)) = \frac{1}{(2+\beta)n-2}.\]

The head of the new edge added to the graph at time $n+1$ is
called the \emph{target vertex} of $v_{n+1}$ and is determined as
follows. If $f_{n+1} = (i,v)$ then the target vertex is $v_i$ and
we say that the choice of target vertex has been made
\emph{uniformly}. If $f_{n+1} = (i,h)$ then the target vertex is
the head of $e_i$ and if $f_{n+1} = (i,t)$ then the target vertex
is the tail of $e_i$, that is $v_i$. When one of the last two
cases occurs, we say that the choice of target vertex has been
made \emph{preferentially} by copying the head or tail, as
appropriate, of $e_i$. Suppose we think of an edge as being
composed of two half-edges so that each half-edge retains one
endpoint of the original edge. Then the target vertex is chosen,
either by choosing one of the $n$ vertices of $G^n_{1,\beta}$
uniformly at random or by choosing one of the $2n-2$ half-edges of
$G^n_{1,\beta}$ uniformly at random and selecting the vertex to
which the half-edge is attached.

The definition implies that for $1 \leq i \leq n$, the probability
that the target vertex of $v_{n+1}$ is $v_i$ is equal to
\begin{equation}\label{eq:fdegree}
\frac{d_n(v_i) + \beta}{(2+\beta)n-2}.\end{equation} We might have
defined $f_{n+1}$ to be a random variable denoting the index of the target
vertex of $v_{n+1}$ and taking probabilities as given in~\eqref{eq:fdegree}. Indeed for much of the sequel we will
abuse notation and assume that we did define $f_{n+1}$ in this
way. However it is useful to have the finer definition when we
prove the concentration results in Section~\ref{sec:conc}.

We extend this model to a random graph process
$(G^n_{m,\beta})$ for $m>1$ as follows: run the graph process
$(G^t_{1,\beta})$ and form $G^n_{m,\beta}$ by taking
$G^{nm}_{1,\beta}$ and merging the first $m$ vertices to form
$v_1$, the next $m$ vertices to form $v_2$ and so on.

Notice that our definition will not immediately extend to the case $\beta=0$
because when $n=1$, the denominator of the expression in~\eqref{eq:fdegree} is zero and so the process cannot start. One way to get around this problem is to define $G^2_{1,0}$ to be the graph with two vertices joined by a single edge and then let the process carry on from there.
A second possibility used in~\cite{Bollobas02}, is to attach an artificial half-edge to $v_1$ at the beginning. This half-edge remains present all through the process so that the sum of the vertex degrees at time $n$ is $2n-1$ rather than $2n-2$ as in the model we use. However it turns out that the choice of which alternative to use makes no difference to the asymptotic form of the expectation of the clustering coefficient and so the results from~\cite{Bollobas02}
are directly comparable with ours.

In the following we only consider properties of the
underlying undirected graph. However, it is helpful to have the
extra notation and terminology of directed graphs to simplify the
reading of some of the proofs.

\section{Subgraphs of $\mathbf{G^n_{1,\beta}}$}\label{sec:smallsub}
Let $S$ be a labelled directed forest with no
isolated vertices, in which each vertex has either one or no
out-going edge and each directed edge $(v_i,v_j)$ has $i >j$.
Moreover if $v_1$ belongs to $S$ than this vertex has no outgoing
edge. The restrictions on $S$ are precisely those that ensure that
$S$ can occur as a subgraph of the evolving M\'ori tree with
$m=1$. We call such an $S$ a \emph{possible forest}.

In this section we generalise the calculation in~\cite{Bollobas02} to calculate the probability that such a graph $S$ is a subgraph of
$G_{1,\beta}^n$ for $\beta > 0$. We will follow the method and
notation of~\cite{Bollobas02} closely.

We emphasise that we are not computing the probability that
$G_{1,\beta}^n$ contains a subgraph isomorphic to $S$; the labels of the vertices of $S$
must correspond to the vertex labels of $G_{1,\beta}^n$ for $S$ to
be considered to be a subgraph of $G_{1,\beta}^n$.

Denote the vertices of $S$ by $v_{s_1},\ldots,v_{s_k}$, where $s_j < s_{j+1}$ for
$1\le j \le k-1 $.
Furthermore, let
\[
V^{-} = \{v_i\in V(S):\text{ there is a }j>i \text{ such that }
(v_j,v_i) \in E(S) \}
\]
and
\[
 V^{+} = \{v_i\in V(S):\text{ there is a }j<i \text{ such that } (v_i,v_j) \in E(S)  \} .
\]
Let $d_S^{\ind}(v)$
($d_S^{\out}(v)$) denote the in-degree (out-degree) of $v$ in $S$.
In particular, $d_S^{\out}(v)$ is either zero or one.
For $t \ge i$, let $R_t(i)=|\{ j>t:(v_j,v_i) \in
E(S) \}|$.
Observe that $R_i(i)=d_S^{\ind}(v_i)$.
Moreover, let $c_S(i)=\sum_{k =1}^{i-1}R_{i-1}(k)$.
Hence $c_S(i)$ is the number of edges in $E(S)$ from
$\{v_i,\ldots,v_n\}$ to $\{v_1,\ldots,v_{i-1} \}$.

\begin{lemma}\label{Le:subgraph1}
Let $\beta >0$ and $S$ be a possible forest. Then for $t \ge s_k$
the probability that $S$ is subgraph of  $G_{1,\beta}^t$ is given by

\begin{align*}
\Pr(S \subset G^t_{1,\beta}) &=
\frac{\beta}{\beta+d_S^{\ind}(v_1)}\prod_{\substack{1\leq i\leq t:\\v_i \in V^-(S)}
}\frac{\Gamma(1+\beta +d_S^{\ind}(v_i))}{\Gamma(1+\beta)}
\\& \phantom{=\
}   \cdot \prod_{\substack{1<i\leq t:\\v_i \in V^{+}}
}\frac{1}{(2+\beta)(i-1)-2}\prod_{\substack{1<i\leq t:\\v_i  \not\in
V^{+}}}\left(1+\frac{c_S(i)}{(2+\beta)(i-1)-2}\right).
\end{align*}
\end{lemma}

\begin{proof}
The proof is a generalisation of the proof for the analogous
result in the case $\beta=0$ in \cite{Bollobas02} but we include
it for completeness.

Let $S_t$ be the subgraph of $S$ induced by the vertices $\{
v_1,\ldots,v_t \}\cap V(S)$. We need to define the following random variables

\[ X_t =  \prod_{(v_l,v_j) \in E(S_t)}I_{(v_l,v_j) \in E(G_{1,\beta}^t)}\prod_{ i \le t} \frac{\Gamma(d_t(v_i)+\beta+R_t(i))}{\Gamma(d_t(v_i)+\beta)}\]
and
\[ Y_t  = \prod_{(v_l,v_j) \in E(S_{t+1})}I_{(v_l,v_j) \in
E(G_{1,\beta}^{t+1})}\prod_{ i \le t}
\frac{\Gamma(d_{t+1}(v_i)+\beta+R_{t+1}(i))}{\Gamma(d_{t+1}(v_i)+\beta)},\]
where $I_A$ is the indicator of the event $A$.

Note  that $d_t(v_j)$ for $1 \le j \le t$ and $X_t$ are functions
of the random variables $f_2,\ldots,f_t$ while $Y_t$ is a function of
the random variables $f_2,\ldots,f_{t+1}$. However, for all $j$, $R_t(j)$ is
deterministic.

Observe that
\[ X_{t+1} = \frac{\Gamma(d_{t+1}(v_{t+1})+\beta +
R_{t+1}(t+1))}{\Gamma(d_{t+1}(v_{t+1})+\beta)} Y_t =
\frac{\Gamma(1+\beta + R_{t+1}(t+1))}{\Gamma(1+\beta)} Y_t.\]

First, assume that there is no $r \leq t$ such that $(v_{t+1},v_r)
\in E(S)$ and so the new edge added at time $t+1$ cannot belong
to $S$. This implies that for $i \leq t$, $R_t(i)= R_{t+1}(i)$ and
\[\prod_{(v_l,v_j) \in E(S_t)}I_{(v_l,v_j)\in E(G^t_{1,\beta})} =
\prod_{(v_l,v_j) \in E(S_{t+1})}I_{(v_l,v_j)\in
E(G^{t+1}_{1,\beta})}.\] Furthermore for all $i \leq t$ with $i
\ne f_{t+1}$, we have $d_{t+1}(v_i) = d_t(v_i)$. We also have
$d_{t+1}(v_{f_{t+1}}) = d_t(v_{f_{t+1}})+1$.

For the moment fix $f_2,\ldots,f_t$ so that $X_t$ is completely
determined. Now,
\[
Y_t =
\left(1+\frac{R_t(f_{t+1})}{d_t(v_{f_{t+1}})+\beta}\right)X_t.
\]
Thus
\begin{align*}
\E{Y_t-X_t|f_2,\ldots,f_t} &=  \sum_{r=1}^t\frac{R_t(r)}{d_t(v_r)+\beta}\Pr(f_{t+1}=r)X_t\\
& = \frac{\sum_{r=1}^tR_t(r)}{(2+\beta)t-2}X_t.
\end{align*}

By taking expectation with respect to $f_2,\ldots,f_t$ we obtain
\[
\E{Y_t} = \left(1+\frac{\sum_{r=1}^t
R_t(r)}{(2+\beta)t-2}\right)\E{X_t}=
\left(1+\frac{c_S(t+1)}{(2+\beta)t-2}\right)\E{X_t}
\]
and
\begin{equation}
\E{X_{t+1}} = \frac{\Gamma(1+\beta+R_{t+1}(t+1))}{\Gamma(1+\beta)}
\left(1+\frac{c_S(t+1)}{(2+\beta)t-2}\right)\E{X_t}. \label{form1}
\end{equation}

Now suppose $(v_{t+1},v_r)$ is an edge of $S$ for some $r < t+1$.
If $f_{t+1} \not= r$ then $X_{t+1}=0$ so we will suppose that
$f_{t+1}=r$. Then for all $i \le t$ with $i \not= r$,
$d_{t+1}(v_i)=d_t(v_i)$, and $d_{t+1}(v_r)=d_t(v_r)+1$.
Furthermore for all $i \le t, i \not= r$ $R_{t+1}(i)=R_t(i)$, but
$R_{t+1}(r)=R_t(r)-1$.

Hence providing $f_{t+1}=v_r$, we have
\[
Y_t=\frac{1}{d_t(v_r)+\beta}X_t.
\]
So
\[
\E{Y_t|f_2,\ldots,f_t}=
\frac{d_t(v_r)+\beta}{(2+\beta)t-2}\frac{X_t}{d_t(v_r)+\beta}=\frac{X_t}{(2+\beta)t-2}.
\]
Thus
\[
\E{X_{t+1}|f_2,\ldots,f_t} = \frac{1}{(2+\beta)t-2}\frac{\Gamma(1+\beta+R_{t+1}(t+1))}{\Gamma(1+\beta)}X_t.\\
\]
So by taking expectation with respect to $f_2,\ldots,f_t$,
\begin{equation}
\E{X_{t+1}} =
\frac{1}{(2+\beta)t-2}\frac{\Gamma(1+\beta+R_{t+1}(t+1))}{\Gamma(1+\beta)}\mathbf{E}[X_t].\label{form2}
\end{equation}

Note that $X_1=
\frac{\Gamma\left(\beta+R_1(1)\right)}{\Gamma\left(\beta \right)}$
and that for $t \ge s_k$, we have $\Pr(S \subset G_{1,\beta}^t)=
\E{X_t}$. Using~\eqref{form1} and~\eqref{form2} and noting that
$R_i(i)=0$ for $v_i \not\in V^{-}$, we have for $t \ge s_k$
\begin{align*}
\Pr(S \subset G^t_{1,\beta}) &=
\frac{\Gamma(\beta+R_1(1))}{\Gamma(\beta)}\prod_{\substack{1<i\leq t:\\v_i \in
V^-}}\frac{\Gamma(1+\beta +R_i(i))}{\Gamma(1+\beta)}
\\& \phantom{=\
} \nonumber\cdot \prod_{\substack{1<i\leq t:\\v_i \in V^{+}}
}\frac{1}{(2+\beta)(i-1)-2}\prod_{\substack{1<i\leq t:\\v_i\not\in
V^{+}}}\left(1+\frac{c_S(i)}{(2+\beta)(i-1)-2}\right).
\end{align*}
This is easily seen to be equivalent to the expression in the
statement of the lemma.
\end{proof}

We now provide a more convenient form for the probability given in
Lemma \ref{Le:subgraph1}. This calculation is almost identical to
the analogous one in \cite{Bollobas02} so we omit the proof.
\begin{lemma}\label{Le:subgraph2}
Let $\beta>0$ and $S$ be a possible forest. Then for $t \geq s_k$
the probability that $S$ is a subgraph of $G^t_{1,\beta}$ is given
by
\begin{align*}
\PullBack{\Pr(S \subset G^{t}_{1,\beta})}
\\ &=
\frac{\beta}{d_S^{\ind}(v_1)+\beta} \prod_{i:v_i \in V^-}
\frac{\Gamma(1+d_S^{\ind}(v_i)+\beta)}{\Gamma(1+\beta)}\\& \phantom{=\
}\cdot \prod_{(v_i,v_j)\in E(S): i> j}
\frac{1}{(2+\beta)(i^{1+\beta}j)^{1/(2+\beta)}}
\exp\left(O\left(\sum_{j=2}^k c_S(s_j)^2/(j-1)\right)\right).
\end{align*}
\end{lemma}
\section{Calculation of Expectations}\label{sec:tri}
Recall that the clustering coefficient $C(G)$ of a graph $G$ is given by
\[ C(G) = \frac{3 \times \text{ number of triangles in }G}{\sum_{v\in V(G)} \binom{d(v)}{2}}.\]
In this section we calculate the expectations of the numerator and
denominator of this expression.

\subsection{Expected Number of Triangles}
We adapt the methods used in \cite{Bollobas02} to the case $\beta
>0$. For fixed $a<b<c$, we first calculate the expected number of
triangles in $G_{m,\beta}^n$ on vertices $v_a,v_b,v_c$. Let
$G_{1,\beta}^{mn}$ be the underlying tree used to form
$G_{m,\beta}^n$. Label the vertices of the tree
$v'_1,\ldots,v'_{mn}$. A triangle on $v_a,v_b,v_c$ arises if there
are vertices $v'_{a_1},v'_{a_2}$ with $(a-1)m+1\leq a_1,a_2\leq
am$, $v'_{b_1},v'_{b_2}$ with $(b-1)m+1\leq b_1,b_2\leq bm$ and
$v'_{c_1},v'_{c_2}$ with $(c-1)m+1\leq c_1,c_2\leq cm$ such that
$v'_{b_1}$ sends its outgoing edge to $v'_{a_1}$, $v'_{c_1}$ sends
its outgoing edge to $v'_{a_2}$ and $v'_{c_2}$ sends its outgoing
edge to $v'_{b_2}$. For this to be possible, we need $c_1\ne c_2$.
Let $S$ be the graph with vertices $v'_{a_1},v'_{a_2},v'_{b_1},v'_{b_2},v'_{c_1},v'_{c_2}$ and edges $(v'_{b_1},v'_{a_1})$, $(v'_{c_1},v'_{a_2})$ and
$(v'_{c_2},v'_{b_2})$. Write $a_1=am-l_1$, $a_2=am-l_2$,
$b_1=bm-l_3$, $b_2=bm-l_4$, $c_1=cm-l_5$ and $c_2=cm-l_6$. The
cases where $a_1=a_2$ and $a_1\ne a_2$ are slightly different. We
concentrate on the former to begin with.

We have $d_S^{\ind}(v_{a_1})=2$, $d_S^{\ind}(v_{b_2})=1$ and otherwise
$d_S^{\ind}(v)=0$. Suppose that $a_1>1$. Then applying
Lemma~\ref{Le:subgraph2} we see
that \begin{equation}\begin{split}\PullBack{\Pr(S\subseteq G_{1,\beta}^{mn})} \\
&=
\frac{\Gamma(3+\beta)\Gamma(2+\beta)}{(\Gamma(1+\beta))^2}\frac{1}{(2+\beta)^3}
\left(\frac{1}{a_1a_2b_2(b_1c_1c_2)^{1+\beta}}\right)^{1/(2+\beta)}\exp(O(1/a)).
\label{eq:tri1}\end{split}\end{equation}
The same expression holds when $a_1=1$ because the extra multiplicative term of $\beta/(2+\beta)$ may be absorbed into the error term.
Note that for $-1\leq x \leq1$, we have $e^x =
1+O(x)$. Furthermore $1/a_i=1/(am)(1+O(1/a))$,
$1/b_i=1/(bm)(1+O(1/a))$ and $1/c_i=1/(cm)(1+O(1/a))$. So we may
rewrite~\eqref{eq:tri1} as follows:
\[
\Pr(S\subseteq G_{1,\beta}^{mn})=
\frac{(1+\beta)^2}{(2+\beta)^2}\frac 1{m^3}
\left(\frac{1}{a^2b^{2+\beta}c^{2+2\beta}}\right)^{1/(2+\beta)}(1+O(1/a)).
\]
In this case where $a_1=a_2$, there are $m^4(m-1)$ ways to choose $a_1,a_2,b_1,b_2,c_1,c_2$ so that there is a corresponding triangle on $v_a,v_b,v_c$ in $G^n_{m,\beta}$.

Now we suppose that $a_1\ne a_2$. We have
$d_S^{\ind}(v_{a_1})=d_S^{\ind}(v_{a_2})=d_S^{\ind}(v_{b_2})=1$ and otherwise
$d_S^{\ind}(v)=0$. Applying Lemma~\ref{Le:subgraph2} and carrying out
similar calculations to those above we obtain
\[
\Pr(S\subseteq G_{1,\beta}^{mn})=
\frac{(1+\beta)^3}{(2+\beta)^3}\frac
1{m^3}\left(\frac{1}{a^2b^{2+\beta}c^{2+2\beta}}\right)^{1/(2+\beta)}(1+O(1/a)).
\]
In this case there are $m^4(m-1)^2$ ways to choose $a_1,a_2,b_1,b_2,c_1,c_2$.

Let $N_{a,b,c}$ denote the number of triangles on
$v_a,v_b,v_c$ in $G_{m,\beta}^n$. From the calculations above, we
see that
\begin{equation}\label{eq:tottri}
\begin{split}
\E{N_{a,b,c}} &= \left(m(m-1)\frac {(1+\beta)^2}{(2+\beta)^2} +
m(m-1)^2\frac {(1+\beta)^3}{(2+\beta)^3}\right)
\left(\frac{1}{a^2b^{2+\beta}c^{2+2\beta}}\right)^{1/(2+\beta)}\\
&\phantom{=\ }\cdot(1+O(1/a)).
\end{split}
\end{equation}

Now let $N$ be the number of triangles in $G_{m,\beta}^n$. Then to
calculate $\E{N}$ we merely sum \eqref{eq:tottri} over all $a,b,c$
with $a<b<c$. If we estimate this sum by integrating, we obtain
the following.
\begin{proposition}\label{pr:extri}
For $\beta>0$, the expected number of triangles in $G^n_{m,\beta}$ is
\[
\left(m(m-1)\frac {(1+\beta)^2}{\beta^2} + m(m-1)^2\frac
{(1+\beta)^3}{\beta^2(2+\beta)}\right) \log n + O(1). \]
\end{proposition}
This result is very different from that obtained
in~\cite{Bollobas02} where it is shown that when $\beta =0$ the
expected number of triangles is $\Theta((\log n)^3)$.

\subsection{Expectation of $\mathbf{\sum_{v\in V(G)}\binom{d(v)}2}$}
We begin by noting that if we regard each edge in the graph as consisting of two half-edges, with each half-edge retaining one endpoint of an edge then $\sum_{v\in V(G_{m,\beta}^n)}\binom{d_n(v)}2$ is the number of pairs of half-edges with the same endpoint. We say such a pair of half-edges is \emph{adjacent}. Suppose that $e_1$ and $e_2$ are half-edges with endpoint $v$. If $e_1$ and $e_2$ form respectively half of edges $vu$ and $vw$ with $u,v,w$ pairwise distinct then we say that $e_1$ and $e_2$ form a \emph{non-degenerate} pair of adjacent half-edges. Otherwise we say that they are \emph{degenerate}.

Calculating the expected number of pairs of adjacent half-edges is
slightly more complicated than calculating the expected number of
triangles because there is less symmetry. We begin by counting the
number of non-degenerate pairs of adjacent half-edges. Let $a<b<c$.
We first
calculate the expected number of pairs $(v_b,v_a)$, $(v_c,v_a)$ of
adjacent half-edges in $G_{m,\beta}^n$ for $\beta
>0$. Just as in the previous section, there are two cases to
consider, and similar calculations, using
Lemma~\ref{Le:subgraph2}, to those above show that the number of
such pairs of adjacent half-edges is
\[
\left(m\frac {1+\beta}{2+\beta} + m(m-1)\frac
{(1+\beta)^2}{(2+\beta)^2}\right)
\left(\frac{1}{a^2b^{1+\beta}c^{1+\beta}}\right)^{1/(2+\beta)}(1+O(1/a)).
\]
By integrating, we see that the total number of pairs of adjacent half-edges in $G_{m,\beta}^n$ for which the common vertex has the smallest index
is
\[
\left(m\frac{2+\beta}{\beta}+m(m-1)\frac{1+\beta}{\beta}\right)n +
O(n^{2/(2+\beta)}).
\]
Now the expected number of pairs $(v_b,v_a)$, $(v_c,v_b)$ of
adjacent half-edges is
\[m^2\frac
{(1+\beta)^2}{(2+\beta)^2}
\left(\frac{1}{ab^{2+\beta}c^{1+\beta}}\right)^{1/(2+\beta)}(1+O(1/a)).\]
Again we integrate to derive that the total number of pairs
of adjacent half-edges in $G_{m,\beta}^n$
for which the common vertex has the middle index
 is $m^2n + O(n^{2/(2+\beta)})$. This is not
surprising because it can be shown that very few vertices either
have loops or do not have $m$ distinct out-neighbours. Each
loopless vertex with $m$ distinct loopless out-neighbours, that
each have $m$ distinct out-neighbours, is the vertex with greatest
index in $m^2$ pairs of adjacent half-edges of this form.

Finally the expected number of pairs $(v_c,v_a)$, $(v_c,v_b)$ of
adjacent half-edges is
\[m(m-1)\frac
{(1+\beta)^2}{(2+\beta)^2}
\left(\frac{1}{abc^{2+2\beta}}\right)^{1/(2+\beta)}(1+O(1/a)).\]
So the total number of pairs of adjacent half-edges in
$G_{m,\beta}^n$ for which the common vertex has the largest index
is $m(m-1)/2n +O(n^{1/(2+\beta)})$. Again this is not surprising
because each loopless vertex with $m$ distinct out-neighbours is
the vertex of greatest index in $\binom{m}{2}$ pairs of adjacent
half-edges of this form.

By carrying out similar calculations to those above, it can be
shown that the number of degenerate pairs of adjacent half-edges
is $O(n^{1/(2+\beta)})$.

Summing over all the possibilities we obtain the following result.
\begin{proposition}\label{pro:exadj}
For $\beta >0$, the expectation of $\sum_{v\in
V(G)}\binom{d(v)}2$ in $G^n_{m,\beta}$ is
\[
\left(\frac{2+5\beta}{2\beta}m^2+\frac{2-\beta}{2\beta}m\right)n +
O(n^{2/(2+\beta)}).\]
\end{proposition}
Again the result is different from that obtained
in~\cite{Bollobas02} where it was shown that for the case $\beta
=0$, the expected number of pairs of adjacent edges is $\Theta(n\log n)$.

\section{Concentration of $\mathbf{\sum_{v\in V(G)}\binom{d(v)}2}$}\label{sec:conc}
In this section we show that the number of pairs of adjacent half-edges
in $G_{m,\beta}^n$ is concentrated about its mean. This justifies
obtaining the clustering coefficient by taking three times the
quotient of the expected number of triangles and the expected
number of pairs of adjacent half-edges. The main strategy is to apply a
variant of the Azuma-Hoeffding inequality from~\cite{Habib98}, by
making use of M\'ori's results~\cite{Mori05} on the evolution of
the maximum degree of $G_{m,\beta}^n$. A key notion in the proof
is to consider the mechanism by which edges incident with a fixed
vertex are added.

Fix $\beta$ and $m$. Let $(H_t)$ be the graph process defined as
follows. Run $(G_{1,\beta}^t)$ and take $H_n$ to be the graph
formed from $G_{1,\beta}^n$ by merging groups of $m$ consecutive
vertices together until there are at most $m$ left and finally
merging the remaining unmerged vertices together. Note that $H_n$
has $\lceil n/m\rceil$ vertices, which we denote by
$v_1,\ldots,v_{\lceil n/m\rceil}$ in the obvious way, and $n-1$
edges. Furthermore, if $m|n$ and the graphs $H_n$ and
$G_{m,\beta}^{n/m}$ are formed from the same instance of the
process $(G_{1,\beta}^t)$, then $H_n$ and $G_{m,\beta}^{n/m}$ are
the same graph.

Let $v_k$ be a vertex of $H_s$ such that $km\leq s$. For $t\geq
s$, we define a partition $\Pi_{k,s}(t)$ of the half-edges
incident with $v_k$. The partition always has $d_s(v_k)+1$ blocks.
When $t=s$, each block of the partition except for one contains one of the
$d_s(v_k)$ half-edges incident with $v_k$; with a slight abuse of
nomenclature the other block, which we call the \emph{base} block,
is initially empty. It follows that if $v_k$ has a loop at time
$s$ then the two half-edges forming the loop are in separate
blocks of $\Pi_{k,s}(s)$. As $t$ increases and more edges are
added to $H$, any newly added half-edge incident with $v_k$ is
added to the partition. If at time $t>s$ the target vertex of the
newly added edge is not $v_k$ then $\Pi_{k,s}(t)=\Pi_{k,s}(t-1)$.
Suppose that at time $t>s$ the target vertex of the newly added
edge $f$ is $v_k$: if $v_k$ is chosen preferentially by copying
the half-edge $e\in A$, where $A$ is a block of $\Pi_{k,s}(t-1)$,
then we form $\Pi_{k,s}(t)$ from $\Pi_{k,s}(t-1)$ by adding the half-edge of $f$
incident with $v_k$ to $A$; if $v_k$ is chosen uniformly then the
half-edge of $f$ incident with $v_k$ is added to the base block.

Suppose that $v_l$ is a vertex of $H_s$ distinct from $v_k$ such
that $lm \leq s$. Suppose further that we choose two distinct
blocks from $\Pi_{k,s}(t)$ and $\Pi_{l,s}(t)$, such that neither
is a base block. The joint distribution of the sizes of the two
blocks is the same for any choice of blocks, whether they are both
chosen from $\Pi_{k,s}(t)$, $\Pi_{l,s}(t)$ or one from each.
Furthermore if we choose either base block from $\Pi_{k,s}(t)$ or
$\Pi_{l,s}(t)$ and one other block that is not a base block, then
again the joint distribution of the sizes of the blocks does not
depend on our choice.

\begin{lemma}\label{le:partition}
Let $v_j$ and $v_k$ be distinct vertices of $H_s$ such that
$\max\{jm,km\}\leq s$. Let $A$ ($B$) be
respectively a block of $\Pi_{j,s}(t)$ ($\Pi_{k,s}(t)$)
such that neither is a base block. Then
\[\E{|A||B|} \leq \E{|A|}\E{|B|} \leq
(t/s)^{2/(2+\beta)}(1+O(1/s)).\]
\end{lemma}
\begin{proof}

Let $e_1,\,e_2$ be half-edges so that at time $s$, $e_1$ is
incident with $v_k$ and $e_2$ is incident with $v_l$. Then let
$a_t$ denote the size, at time $t$, of the block of $\Pi_{k,s}(t)$
containing $e_1$ and let $b_t$ be defined similarly with respect
to $\Pi_{l,s}(t)$ and $e_2$. We first establish the second
inequality. We have $\E{a_s}=1$ and for $t\geq s$,
\begin{align}\label{eq:expect}
 \E{a_{t+1}|a_t} = a_t\left(1 + \frac{1}{(2+\beta)t-2}\right).
\end{align}
Hence
\[ \E{a_{t+1}} =  \frac{t-1/(2+\beta)}{t-2/(2+\beta)}\E{a_t}.\]
Solving this recurrence, we obtain
\[\E{a_t} =  \frac{\Gamma\left(t-\frac 1{2+\beta}\right)\Gamma\left(s-\frac 2{2+\beta}\right)}{\Gamma\left(t-\frac 2{2+\beta}\right)\Gamma\left(s-\frac
1{2+\beta}\right)}.\] A standard result on the ratio of gamma
functions \cite{Marichev82} states that if $a,b$ are fixed members
of $\reals$ then for all $x>\max\{|a|,|b|\}$,
\[
\frac{\Gamma(x+b)}{\Gamma(x+a)} = x^{b-a}(1+O(1/x)) .\] Using this
result, we obtain
\[ \E{a_t} \leq (t/s)^{1/(2+\beta)} (1+O(1/s)).\]
Since $|A|$ and $|B|$ are identically distributed, the second
inequality in the lemma follows. We prove the first inequality by
using induction on $t$. Observe that $(a_{t+1},b_{t+1})$ can take
the values $(a_t+1,b_t)$, $(a_t,b_t+1)$ and $(a_t,b_t)$ with
probabilities respectively $a_t/((2+\beta)t-2)$,
$b_t/((2+\beta)t-2)$ and $1-(a_t+b_t)/((2+\beta)t-2)$. Therefore
\[
\E{a_{t+1}b_{t+1}|a_tb_t} = a_tb_t + \frac{2a_tb_t}{(2+\beta)t-2}\]
and from \eqref{eq:expect} we get
\[\E{a_{t+1}}\E{b_{t+1}} = \E{a_t}\E{b_t}\left(1 +
\frac{1}{(2+\beta)t-2}\right)^2.
\]
So
\[\E{a_{t+1}b_{t+1}}-\E{a_{t+1}}\E{b_{t+1}} \leq \left(1+\frac{2}{(2+\beta)t-2}\right)(\E{a_tb_t}-\E{a_t}\E{b_t})\]
and hence the result follows by induction.
\end{proof}

When the maximum degree of $H_t$ becomes unusually large and the
target vertex is chosen to be a vertex of maximum degree, the
number of pairs of adjacent edges increases by an unusually large
amount. The next result enables us to show that the probability of
this happening is extremely small. Let $\Delta(G)$ denote the
maximum degree of $G$. The following is a very slight
reformulation of what M\'ori proves in~\cite[Theorem~3.1]{Mori05}.
\begin{theorem}
For any positive integer $k$, there exists $\tilde M_k \in
\reals$, such that for all $n$,
\[\mathbf{E}\left[\left(\frac{\Delta(G_{1,\beta}^n)+\beta}{n^{1/(2+\beta)}}\right)^k\right]
\leq \tilde{M}_k.\]
\end{theorem}
The following corollary is straightforward.
\begin{corollary}\label{Co:universalM}
For any positive integers $k,m$, there exists $M_{k,m} \in \reals$
such that for all positive integers $i_1,\ldots,i_k$,
\[ \mathbf E \left[\frac{\Delta(H_{mi_1})}{(mi_1)^{1/(2+\beta)}} \cdots \frac{\Delta(H_{mi_k})}{(mi_k)^{1/(2+\beta)}} \right ] \leq M_{k,m}.\]
\end{corollary}
\begin{proof}
Since $\Delta(H_{mi_1}),\ldots,\Delta(H_{mi_k})$ are all positive
we have
\[\frac{\Delta(H_{mi_1})}{(mi_1)^{1/(2+\beta)}} \cdots \frac{\Delta(H_{mi_k})}{(mi_k)^{1/(2+\beta)}} \leq \sum_{j=1}^k \left(\frac{\Delta(H_{mi_j})}{(mi_j)^{1/(2+\beta)}}\right)^k\]
and so
\[\mathbf{E}\left[\frac{\Delta(H_{mi_1})}{(mi_1)^{1/(2+\beta)}} \cdots \frac{\Delta(H_{mi_k})}{(mi_k)^{1/(2+\beta)}}\right] \leq \sum_{j=1}^k\mathbf{E} \left[\left(\frac{\Delta(H_{mi_j})}{(mi_j)^{1/(2+\beta)}}\right)^k\right].\]
Recall that $H_{mi}$ is formed by merging together blocks of $m$
consecutive vertices in an instance of $G_{1,\beta}^{mi}$. So we
have $\E{(\Delta(H_{mi}))^k} \leq
\E{(m\Delta(G_{1,\beta}^{mi}))^k}$. Hence
\[
\sum_{j=1}^k\mathbf{E}
\left[\left(\frac{\Delta(H_{mi_j})}{(mi_j)^{1/(2+\beta)}}\right)^k\right]
\leq m^k \sum_{j=1}^k \mathbf{E}
\left[\left(\frac{\Delta(G_{1,\beta}^{mi_j})}{(mi_j)^{1/(2+\beta)}}\right)^k\right]
\leq km^k\tilde M_{k}.\] The result follows by taking
$M_{k,m}=km^k\tilde{M}_k$.
\end{proof}

Before we can state the large deviation result that we use, we
need some more definitions. Recall that $f_i$ is a random variable
which determines the index of the target vertex of $v_i$ and that
the values taken by $f_2,f_3,\ldots,f_t$ together determine $H_t$.
Furthermore the set of values that $f_i$ can take is denoted by
$\Omega_i$ and $f_2,\ldots,f_t$ are independent. Let
$\Omega=\prod_{i=2}^t \Omega_i$.

Let $\vect X=(f_2,\ldots,f_t)$. We let $H_t(\vect X)$ be the
instance of $H_t$ determined by the random variables
$f_2,\ldots,f_t$. We will also use this notation both for other
random variables associated with $H_t$ and when some or all of the
$f_i$'s are set to a particular value. The meaning should be clear
from the context but we will generally use $\omega_i$ for a member
of $\Omega_i$ and $f_i$ for a random variable taking values in
$\Omega_i$.

Let $D(\vect X) = \sum_{v \in V(H_t(\mathbf{X}))} \binom{d_t(v)}2$
and let $F(\vect X) = D(\vect X) t^{-2/(2+\beta)}$. Now let
$g:\prod_{i=2}^s \Omega_i \rightarrow \reals$ such that
\[(\omega_2,\ldots,\omega_s) \mapsto \E{F(\omega_2,\ldots,\omega_s,f_{s+1},\ldots,f_t)}\]
and let $\ran:\prod_{i=2}^{s-1} \Omega_i \rightarrow \reals$ such
that
\[ (\omega_2,\ldots,\omega_{s-1}) \mapsto \sup{\{|g(\omega_2,\ldots,\omega_{s-1},x)-g(\omega_2,\ldots,\omega_{s-1},y)|: x,y \in \Omega_s\}}.\]
So $\ran(\omega_2,\ldots,\omega_{s-1})$ measures the maximum
amount that the expected value of $F(\vect X)$ changes when the
value of $f_s$ is changed.

For $\vect \omega \in \Omega$, let
\[ R^2(\vect \omega) = \sum_{k=2}^{t} \ran(\omega_2,\ldots,\omega_{k-1})^2.\]
Our aim is to bound $R^2(\vect \omega)$ as $\vect \omega$ runs
over all members of $\Omega$ with the possible exception of those
belonging to some `bad' subset $\mathcal B$ which we hope to have
small probability. We specify $\mathcal B$ below but for the
moment let $\mathcal B$ be any subset of $\Omega$. Let
\[ r^2 = \sup \{R^2(\vect \omega) : \vect \omega \in \Omega \setminus \mathcal{B}\}.\]
Then Theorem 3.7 in \cite{Habib98} yields the following inequality. For all $x>0$,
\[ \Pr(|F(\vect X) - \E{F(\vect X)}| \geq x) \leq 2(e^{-2x^2/r^2} + \Pr(\vect X \in \mathcal{B})).\]

Fix $\delta>0$. We let
\[ \mathcal{B}_\delta = \left\{\vect X\in \Omega: \sum_{i=1}^n \left(\frac{\Delta(H_{mi}(\vect X))}{(mi)^{2/(2+\beta)}}\right)^2 \geq n^{\frac{\beta}{2+\beta}+\delta}\right\}.\]
Then we have the following.
\begin{lemma}\label{le:maxdegree}
For any $\delta>0$ and $\gamma>0$, there exists $L$ such that $\Pr(\mathcal{B}_\delta) \le
L\frac{1}{n^{\gamma}}$, where $L$ is a constant
depending on $\delta,\gamma,\beta,m$ but not on $n$.
\end{lemma}
\begin{proof}
For any positive integer $k$,
Markov's inequality gives
\[\Pr(\mathcal{B}_\delta) \leq  \frac{\mathbf{E}\left[\left(\sum_{i=1}^n\left(\frac{\Delta(H_{mi}(\vect X) )}{(mi)^{2/(2+\beta)}}\right)^2\right)^k\right]} {n^{\frac{\beta k}{2+\beta}+k\delta}}.
\]
The numerator of this fraction is equal to
\[
\mathbf{E}\left[\sum_{i_1=1}^n\cdots\sum_{i_k=1}^n \left( \frac{\Delta(H_{mi_1}(x))}{(mi_1)^{1/(2+\beta)}}\right)^2\cdots\left( \frac{\Delta(H_{mi_k}(x))}{(mi_k)^{1/(2+\beta)}}\right)^2 \frac{1}{(m^ki_1 \cdots i_k)^{2/(2+\beta)}}\right]. \\
\]
Using Corollary~\ref{Co:universalM} this is at most
\begin{align*}
M_{2k,m} \sum_{i_1=1}^n\cdots\sum_{i_k=1}^n \left( \frac{1}{(m^ki_1
\cdots i_k)^{2/(2+\beta)}}\right) &= M_{2k,m}\left( \sum_{i=1}^n
\frac{1}{(mi)^{\frac{2}{2+\beta}}}\right)^k\\
&\le M_{2k,m}\left( \frac{2+\beta}{\beta}
\frac{n^{\frac{\beta}{2+\beta}}}{m^{\frac{2}{2+\beta}}} \right)^k.
\end{align*}
Hence
\begin{align*}\Pr(\mathcal{B}_\delta) &\leq \frac
{M_{2k,m}\left(\frac{2+\beta}{\beta}\frac{1}{m^{\frac{2}{2+\beta}}}\right)^k}{n^{k\delta}}
\end{align*}
and so letting $k=\lceil \gamma/\delta \rceil$ gives the result.
\end{proof}

We can now state the main result of this section concerning the
concentration of the number of pairs of adjacent half-edges around
its expectation.

\begin{theorem}\label{th:conc}
Let $\beta>0$. For any $\epsilon >0$, the number $D$ of pairs of
adjacent half-edges in $G_{m,\beta}^n$ is concentrated about its
expected value within $O(n^{(4+\beta)/(4+2\beta)+\epsilon})$. More
precisely, for any $\epsilon >0$ and $\gamma>0$ there exists $n^*$
such that for all $n \ge n^*$
\[ \Pr\left(|D- \E{D}| \geq n^{\frac{4+\beta}{4+2\beta}+\epsilon}\right) \le \frac{1}{n^{\gamma}}.\]
\end{theorem}
\begin{proof}
Let $t=nm$, and fix $s\leq t$. Let $s'=m\lceil s/m \rceil$, so we have $s'\leq t$.
Now let
\[\vect \omega_x=(\omega_2,\ldots,\omega_{s-1},x,\omega_{s+1},\ldots,\omega_t) \quad
\text{and} \quad \vect \omega_y =
(\omega_2,\ldots,\omega_{s-1},y,\omega_{s+1},\ldots,\omega_t),\]
where $\omega_i\in\Omega_i$ and $x,y\in \Omega_s$.
For $z\in \{x,y\}$, let $d^z_{t}(v)$ denote the total degree of
$v$ at time $t$ in $H_{t}(\vect \omega_z)$ and let $e$ denote
the edge added at time $s$. Suppose that in $H_t(\vect \omega_x)$
the target vertex of $e$ is $v_{k_1}$ and in $H_t(\vect \omega_y)$
the target vertex of $e$ is $v_{k_2}$. Note that at any time, for
every vertex $v$ other than $v_{k_1}$ or $v_{k_2}$, the degree of
$v$ is the same in $H_t(\vect \omega_x)$ and $H_t(\vect
\omega_y)$. Therefore $F(\vect \omega_x)-F(\vect \omega_y)$
depends only on the degrees of $v_{k_1}$ and $v_{k_2}$ and is
given by
\begin{equation}\label{eq:Fdiff}\begin{split}\PullBack{F(\vect \omega_x)-F(\vect \omega_y)}\\ &= t^{-2/(2+\beta)}
\left( \binom{d^x_t(v_{k_1})}2 + \binom{d^x_t(v_{k_2})}2 -
\binom{d^y_t(v_{k_1})}2  -\binom{d^y_t(v_{k_2})}2
\right).\end{split}\end{equation} From now on we will assume that $k_1 \ne
k_2$, because otherwise $F(\vect \omega_x)-F(\vect \omega_y)=0$.
Consider the changes that occur to $H_{s'}$ if we replace $\vect
\omega_y$ by $\vect \omega_x$. First the head of $e$ is moved from
$v_{k_2}$ to $v_{k_1}$. Second it is possible that each of the at
most $m-1$ edges that are added in the time interval $[s+1,s']$
also have an endpoint moved from $v_{k_2}$ to $v_{k_1}$: this will
happen if the target vertex of an edge added in the interval
$[s+1,s']$ is chosen by preferentially copying the head of an edge
which has been moved from $v_{k_2}$ to $v_{k_1}$, in particular if
the target vertex is chosen by preferentially copying the head of
$e$. Consequently we have
\[ d^y_{s'}(v_{k_1})+1 \leq d_{s'}^x(v_{k_1}) \leq d_{s'}^y(v_{k_1})+m\]
and furthermore
\[ d_{s'}^x(v_{k_1}) + d_{s'}^x(v_{k_2}) = d_{s'}^y(v_{k_1}) + d_{s'}^y(v_{k_2}).\]
Let $d=d^x_{s'}(v_{k_1})-d^y_{s'}(v_{k_1})$,
$d_1=d_{s'}^y(v_{k_1})$ and $d_2=d_{s'}^x(v_{k_2})$. Note that
both $d_1$ and $d_2$ and consequently also $|d_1-d_2|$ are at most
$\Delta(H_{s-1}(\omega_1,\ldots,\omega_{s-1}))+m$.

Now let $A_0,A_1,\ldots,A_{d_1}$, ($B_0,B_1,\ldots,B_{d_2}$)
denote the blocks of the partition  $\Pi_{k_1,s'}(t)$ in
$H_t(\vect \omega_y)$ ($\Pi_{k_2,s'}(t)$ in $H_t(\vect \omega_x)$)
with $A_0$ ($B_0$) denoting the base block. The partition
$\Pi_{k_1,s'}(t)$ in $H_t(\vect \omega_x)$ contains the blocks
$A_0,\ldots,A_{d_1}$ but also $d$ further blocks which we label
$C_1,\ldots,C_d$. Then the partition $\Pi_{k_2,s'}(t)$ in
$H_t(\vect \omega_y)$ contains the blocks
$B_0,\ldots,B_{d_2}$, $C_1,\ldots,C_d$. So using~\eqref{eq:Fdiff}, we
have
\begin{equation}\label{eq:Fdiff1} F(\vect \omega_x)-F(\vect \omega_y)
 = t^{-2/(2+\beta)} \left( \sum_{i=0}^{d_1}
  \sum_{j=1}^d |A_i||C_j|-\sum_{i=0}^{d_2} \sum_{j=1}^d |B_i||C_j|  \right).\end{equation}

Now let
\begin{align*}
\vect \omega_x &=
(\omega_2,\ldots,\omega_{s-1},x,\omega_{s+1},\ldots,\omega_{s'},f_{s'+1},\ldots,f_t)\\
\intertext{and} \vect \omega_y &=
(\omega_2,\ldots,\omega_{s-1},y,\omega_{s+1},\ldots,\omega_{s'},f_{s'+1},\ldots,f_t).\end{align*}
So both $H_t(\vect \omega_x)$ and $H_t(\vect \omega_y)$ evolve
deterministically until time $s'$ but randomly thereafter.

Recall that $d\leq m$ and that $|d_1-d_2|$ is at most
$\Delta(H_{s-1}(\omega_2,\ldots,\omega_{s-1}))+m$. Hence from~\eqref{eq:Fdiff1}, Lemma~\ref{le:partition} and the remarks immediately preceding the
lemma, we see that
\[ |\E{F(\vect \omega_x)-F(\vect \omega_y)}|  \leq
(\Delta(H_{s-1}(\omega_2,\ldots,\omega_{s-1}))+m)m
(1/s')^{2/(2+\beta)}(1+O(1/s')).\] Notice that this expression
does not depend on $x$ or $y$ and holds for every
$\omega_{s+1},\ldots,\omega_{s'}$. Consequently
\[ \ran(\omega_2,\ldots,\omega_{s-1}) \leq
(\Delta(H_{s-1}(\omega_2,\ldots,\omega_{s-1}))+m)m
(1/s')^{2/(2+\beta)}(1+O(1/s')).\] Now let $\vect \omega \in
\Omega \setminus \mathcal{B}_\delta$. Then
\begin{align*}
R^2(\vect \omega) &= \sum_{s=2}^{nm}
(\Delta(H_{s-1}(\omega_2,\ldots,\omega_{s-1}))+m)^2m^2
(1/s')^{4/(2+\beta)}(1+O(1/s'))\\
&\leq m^2\sum_{s=2}^{nm}
\left(\frac{2\Delta(H_{s'}(\omega_2,\ldots,\omega_{s'}))}{s'^{2/2+\beta}}\right)^2(1+O(1/s'))\\
&\leq 4m^3\sum_{i=1}^n
\left(\frac{\Delta(H_{mi}(\omega_2,\ldots,\omega_{mi}))}{(mi)^{2/2+\beta}}\right)^2(1+O(1/i'))\\
&\leq cn^{\frac{\beta}{2+\beta}+\delta},
\end{align*}
where $c$ is a constant.

Hence
\begin{align*}
\Pr\left(|D(\mathbf{X})-\E{D(\mathbf{X})}| \geq
n^{\frac{4+\beta}{4+2\beta}+\epsilon}\right)
&= \Pr\left(|F(\vect X)-\E{F(\vect X)}| \geq n^{\frac{\beta}{4+2\beta}+\epsilon}\right)\\
&\leq 2\exp\left(\frac{-2
n^{\frac{\beta}{2+\beta}+2\epsilon}}{cn^{\frac{\beta}{2+\beta}+\delta}}\right)+2\Pr(\mathcal{B}_\delta).
\end{align*}
If we choose $\delta=\epsilon$ then the first term is at most
$\frac{1}{2n^{\gamma}}$ for any $\gamma >0$ and sufficiently large
$n$. Applying Lemma \ref{le:maxdegree} with any
$\gamma^{*}>\gamma$ we see that for sufficiently large $n$ we also
have $2\Pr(\mathcal{B_{\epsilon}}) \le \frac{1}{2n^{\gamma}}$.
Hence the result follows.
\end{proof}

\section{Expected clustering coefficient}\label{sec:final}
In this section we finally state and prove our main result.
\begin{theorem}
For any $\beta>0$, the expected clustering coefficient of
$G_{m,\beta}^n$ is given by
\[
\mathbf{E}[C(G_{m,\beta}^n)]=
\frac{3c_1\log n}{c_2n}+O(1/n),
\]
where \[c_1=m(m-1)\frac {(1+\beta)^2}{\beta^2} + m(m-1)^2\frac
{(1+\beta)^3}{\beta^2(2+\beta)} \] and \[
c_2=\frac{2+5\beta}{2\beta}m^2+\frac{2-\beta}{2\beta}m.\]
\end{theorem}

\begin{proof}
Recall that $N=N(G_{m,\beta}^n)$, $D=D(G_{m,\beta}^n)$ denote
respectively the number of triangles and pair of adjacent edges in
$G_{m,\beta}^n$. The expected clustering coefficient is given by
$\E{C(G_{m,\beta}^n)} = \E{3N/D}$.

Choose $\epsilon$ so that $0 < \epsilon < \frac{\beta}{4+2\beta}$
and let $\eta=\epsilon+\frac{4+\beta}{4+2\beta}<1$. Let $I$ denote
the interval $[\E{D}-n^{\eta},\E{D}+n^{\eta}]$.
From Proposition~\ref{pro:exadj} we have $\E{D}-n^{\eta} = c_2n - (1+o(1))n^{\eta}$ and
$\E{D}+n^{\eta} = c_2n + (1+o(1))n^{\eta}$. Let $n \geq n^{\ast}$, the minimum value of $n$ such that Theorem~\ref{th:conc} may be applied with $\gamma=4$.
Since $C(G_{m,\beta}^n) \leq m$, an upper
bound for $\E{C(G_{m,\beta}^n)}$ may be obtained as
follows.
\begin{align*}
\E{C(G_{m,\beta}^n)} &\leq \sum_{j=1}^{\infty}\sum_{i \in I
}\frac{3j}{i}\Pr(N=j,D=i)+m\Pr(D \not\in I)\\
&\leq \sum_{j=1}^{\infty}\frac{3j}{c_2n-(1+o(1))n^{\eta}}\Pr(N=j)+m\Pr(D \not\in I).
\end{align*}
Applying Theorem~\ref{th:conc} with $\gamma=1$ and then Proposition~\ref{pr:extri}, we obtain
\begin{align*}
\E{C(G_{m,\beta}^n)} &\le
\sum_{j=1}^{\infty}\frac{3j}{c_2n-(1+o(1))n^{\eta}}\Pr(N=j)+\frac{m}{n}\\
&= \frac{3c_1\log n}{c_2 n}(1+(1/c_2+o(1))n^{\eta-1})+\frac{m}{n}\\
&= \frac{3c_1\log n}{c_2 n} + O(1/n).
\end{align*}

A lower bound for $\mathbf{E}(C(G_{m,\beta}^n))$ may be obtained
as follows.
\begin{align*}
\E{C(G_{m,\beta}^n)} &\ge \sum_{j=1}^{\infty}\sum_{i \in I
}\frac{3j}{i}\Pr(N=j,D=i)\\
&\ge \sum_{j=1}^{\infty}\sum_{i \in I
}\frac{3j}{c_2n + (1+o(1))n^{\eta}}\Pr(N=j,D=i)\\
&=\frac{3\E{N}}{c_2n + (1+o(1))n^{\eta}}\\
&\phantom{=} {}-\sum_{j=1}^{\infty}\sum_{i
\not\in I }\frac{3j}{c_2n + (1+o(1))n^{\eta}}\Pr(N=j,D=i).
\end{align*}
Now since there are at most $n^3m^3$ triangles in $G_{m,\beta}^n$
\begin{align*}
\sum_{j=1}^{\infty} \sum_{i \not\in
I}\frac{3j}{c_2n + (1+o(1))n^{\eta}}\Pr(N=j,D=i) \le
\frac{3n^3m^3}{c_2n + (1+o(1))n^{\eta}}\Pr(D \not\in I).
\end{align*}
Applying Theorem \ref{th:conc} with $\gamma=4$ shows that this is
$O\left(1/n \right)$. Finally
\[
\frac{3\E{N}}{c_2n + (1+o(1))n^{\eta}}= \frac{3c_1 \log n}{c_2
n}(1-(1/c_2+o(1))n^{\eta-1}) = \frac{3c_1 \log n}{c_2 n} +
O(1/n).\]
\end{proof}
\section{Conclusion}
Our main result shows that for $\beta >0$ the expectation of the clustering coefficient of the M\'ori graph is asymptotically proportional to $\log n/n$ and consequently that the M\'ori graphs do not have the small-worlds property.
Bollob\'as and Riordan showed for an almost identical model that when $\beta=0$, the expectation of the clustering coefficient is asymptotically proportional to $(\log n)^2/n$.
An unexpected consequence, for which we do not yet have a good explanation, is that the clustering coefficient has a discontinuity at $\beta =0$.

\begin{section}{Acknowledgement}
Some of this research was carried out while the first author was
visiting the University of Bordeaux. We thank Philippe Duchon and
Nicolas Hanusse for their hospitality and for useful discussions.
 We also thank John Harris for his
patient explanation of martingales.
\end{section}


\begin{thebibliography}{99}

\bibitem{Barabasi99}
A.-L. Barab\'asi and R.~Albert.
\newblock Emergence of scaling in random networks.
\newblock {\em Science}, 286(5439):509--512, 1999.

\bibitem{Albert02}
A.-L. Barab\'asi and R.~Albert.
\newblock Statistical mechanics of complex networks.
\newblock {\em Reviews of modern Physics}, 74:47--97, 2002.

\bibitem{Bollobas02}
B.~Bollob\'as and O.~M. Riordan.
\newblock Mathematical results on scale-free random graphs.
\newblock In S.~Bornholdt and H.~G. Schuster, editors, {\em Handbook of Graphs
  and Networks: From the Genome to the Internet}, chapter~1, pages 1--34.
  Wiley-VCH, Berlin, 2003.

\bibitem{Bollobas01}
B.~Bollob\'as, O.~M. Riordan, J.~Spencer, and G.~Tusn\'ady.
\newblock The degree sequence of a scale-free random graph process.
\newblock {\em Random Structures and Algorithms}, 18(3):279--290, 2001.

\bibitem{Buckley04}
P.~G. Buckley and D.~Osthus.
\newblock Popularity based random graph models leading to a scale-free degree
  sequence.
\newblock {\em Discrete Mathematics}, 282:53--63, 2004.

\bibitem{Cooper03}
C.~Cooper and A.~Frieze.
\newblock A general model of web graphs.
\newblock {\em Random Structures and Algorithms}, 22(3):311--335, 2003.

\bibitem{Durret06}
R.~Durrett.
\newblock {\em Random Graph dynamics}.
\newblock Cambridge University Press, 2006.

\bibitem{Marichev82}
O.~I. Marichev.
\newblock {\em Handbook of integral transforms of higher transcendental
  functions, theory and algorithmic tables}, chapter The gamma function and its
  properties, pages 43--52.
\newblock Ellis Horwood Limited, 1983.

\bibitem{Habib98}
C.~McDiarmid.
\newblock Concentration.
\newblock In {\em Probabilistic Methods for Algorithmic Discrete Mathematics},
  number~16 in Algorithms and Combinatorics, pages 195--248. Springer, 1998.

\bibitem{Mori03}
T.~F. M\'ori.
\newblock On random trees.
\newblock {\em Studia Scientiarum Mathematicarum Hungarica}, 39(1-2):143--155,
  2002.

\bibitem{Mori05}
T.~F. M\'ori.
\newblock The maximum degree of the {B}arab\'asi-{A}lbert random tree.
\newblock {\em Combinatorics, Probability and Computing}, 14:339--348, 2005.

\bibitem{Watts98}
D.~J. Watts and S.~H. Strogatz.
\newblock Collective dynamics of `small-world' networks.
\newblock {\em Nature}, 393(6684):440--442, 1998.

\end{thebibliography}
\end{document}